\documentclass[12pt]{amsart}

\usepackage{amssymb}
\usepackage[dvips]{graphicx}
\usepackage{wrapfig}

\newtheorem{thm}{Theorem}[section]
\newtheorem{cor}[thm]{Corollary}
\newtheorem{lem}[thm]{Lemma}
\newtheorem{prop}[thm]{Proposition}

\theoremstyle{definition}

\newtheorem{eg}[thm]{Example}
\theoremstyle{remark}
\newtheorem{rem}[thm]{Remark.}

\title{Region crossing change is an unknotting operation}
\author{Ayaka SHIMIZU}
\address{Department of Mathematics, Hiroshima University, 1-3-1 Kagamiyama, Higashi-Hiroshima, 7398526, JAPAN}
\email{shimizu1984@gmail.com}
\date{\today}

\begin{document}

\maketitle

\begin{abstract}
A region crossing change is a local transformation on a knot or link diagram. 
We show that a region crossing change on a knot diagram is an unknotting operation, 
and we define the region unknotting number for a knot diagram and a knot. 
\end{abstract}

\section{Introduction}

\begin{figure}[ht]
\begin{center}
\includegraphics[width=135mm]{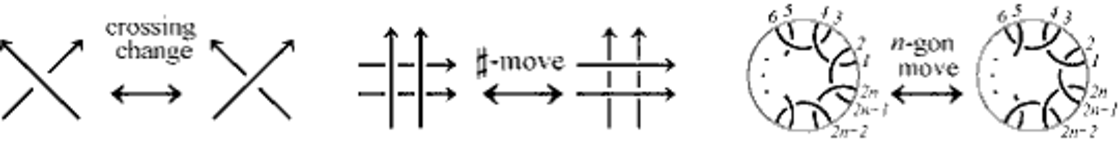}
\caption{}
\label{local}
\end{center}
\end{figure}

\begin{figure}[ht]
\begin{center}
\includegraphics[width=50mm]{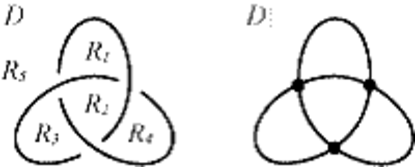}
\caption{}
\label{rcc-ex2}
\end{center}
\end{figure}

\begin{figure}[ht]
\begin{center}
\includegraphics[width=110mm]{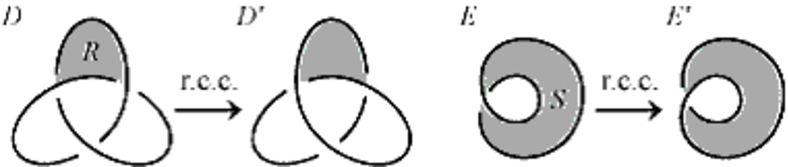}
\caption{}
\label{rcc-ex}
\end{center}
\end{figure}

An \textit{unknotting operation} is a local transformation of a knot diagram such that any diagram can be transformed into a diagram of the trivial knot by a finite sequence of these operations. 
Unknotting operations play an important role in knot theory, and many unknotting operations have been studied. 
For example, it is well-known that the crossing change, indicated in Figure \ref{local}, is an unknotting operation. 
It is also known that the \textit{$\sharp$-move}, indicated in Figure \ref{local}, is an  unknotting operation \cite{murakami}, and 
that an \textit{$n$-gon move}, indicated in Figure \ref{local}, is an unknotting operation \cite{miyazawa}, \cite{nakanishi}. 
Let $D$ be a link diagram on $S^2$, and $\vert D\vert$ be the four-valent graph obtained from $D$ by replacing each crossing with a vertex. 
We call each component of $S^2-\vert D\vert$ a \textit{region} of $D$. 
A diagram $D$ with $c$ crossings has $2c$ edges and, therefore, $c+2$ regions because the Euler characteristic of $S^2$ is $2$. 
For example, the diagram $D$ with three crossings in Figure \ref{rcc-ex2} has six edges and five regions $R_1,R_2,\dots$ and $R_5$. 
A \textit{region crossing change} at a region $R$ of $D$ is the local transformation on $D$ by the changing 
all the crossings on $\partial R$. 
For example in Figure \ref{rcc-ex}, we obtain $D'$ (resp. $E'$) from $D$ (resp. $E$) by applying a region crossing change at $R$ (resp. $S$). 
A $\sharp$-move and an $n$-gon move on a knot diagram are examples of region crossing changes. 
We remark that we can apply a region crossing change on a non-alternating region even though we cannot apply an $n$-gon move. 
The region crossing change was proposed by K. Kishimoto in a seminar at Osaka City University in 2010. 
Kishimoto raised the following question: 
\textit{Is a region crossing change on a knot diagram an unknotting operation?} 
This means, ``Can we transform any diagram into a diagram of the trivial knot by region crossing changes''. 
We will prove the following theorem in this paper: 

\phantom{x}
\begin{thm}
Let $D$ be a knot diagram, and $c$ a crossing point of $D$. 
Let $D'$ be the diagram obtained from $D$ by the crossing change at $c$. 
Then, there exist region crossing changes which transform $D$ into $D'$. 
\label{mainthm}
\end{thm}
\phantom{x}

\noindent The proof is given in Section \ref{proof-s}. 
Since a crossing change on a knot diagram is an unknotting operation, we have the following corollary of Theorem \ref{mainthm} which answers Kishimoto's question: 

\phantom{x}
\begin{cor}
A region crossing change on a knot diagram is an unknotting operation. 
Therefore, we can transform any diagram into a diagram of the trivial knot by region crossing changes. 
\label{answer}
\end{cor}
\phantom{x}

\begin{rem}
For a link diagram, the answer to Kishimoto's question is negative. 
For example, the link diagram in Figure \ref{hopf} can not be transformed into a diagram of a trivial link by any number of region crossing changes. 
\begin{figure}[ht]
\begin{center}
\includegraphics[width=25mm]{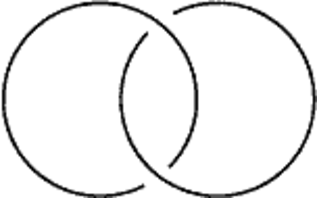}
\caption{}
\label{hopf}
\end{center}
\end{figure}
\end{rem}

\noindent We define the \textit{region unknotting number} $u_R(D)$ of a knot diagram $D$ to be the minimum number of 
region crossing changes necessary to obtain a diagram of the trivial knot from $D$. 
For example, we have $u_R(D)=1$ for the diagram $D$ in Figure \ref{rcc-ex}. 
H. A. Miyazawa showed in \cite{miyazawa} that for any knot $K$, there exists an integer $n$ such that a diagram of $K$ can be 
transformed into a diagram of the trivial knot by one $n$-gon move. 
Therefore, every knot $K$ has a diagram $D$ such that $u_R(D)=1$. 
We define the \textit{region unknotting number} $u_R(K)$ of a knot $K$ to be the minimal $u_R(D)$ for all minimal crossing diagrams $D$ of $K$. 
We have the following theorem: 

\phantom{x}
\begin{thm}
Let $K$ be a knot and $c(K)$ be the crossing number of $K$, then 
$$u_R(K)\le \frac{c(K)}{2}+1.$$
\label{urk-thm}
\end{thm}
\phantom{x}

\noindent The proof is given in Section \ref{run}. 
The rest of this paper is organized as follows: 
In Section \ref{lem-s}, we develop the properties of region crossing changes that are used in proving Theorem \ref{mainthm}. 
In Section \ref{proof-s}, we prove Theorem \ref{mainthm}. 
In Section \ref{run}, we consider the region unknotting number and prove Theorem \ref{urk-thm}. 
In the appendix, we discuss minimal crossing diagrams on $S^2$ for prime alternating knots.

\section{Properties of region crossing changes}
\label{lem-s}

In this section, we discuss the properties of region crossing changes on a link diagram. 
Let $D$ be a link diagram and $R$ a region of $D$. 
We denote by $D(R)$ the diagram obtained from $D$ by the region crossing change on $R$. 
For two regions $R_1$ and $R_2$ of $D$, we denote by $D(R_1, R_2)$ the diagram obtained from $D$ by the region crossing changes first on $R_1$, and then on $R_2$. 
We have $D(R_1, R_2)=D(R_2, R_1)$ and $D(R, R)=D$ because the result of crossing changes does not depend on the order, 
and two crossing changes at a crossing point cancel. 
For regions $R_1,R_2,\dots$ and $R_n$ of $D$, the set of regions $P=R_1\cup R_2\cup \dots \cup R_n$ 
allows us to denote by $D(P)$ the diagram obtained from $D$ by region crossing changes on $R_1,R_2,\dots$ and $R_n$. 
We have the following lemma: 

\phantom{x}
\begin{lem}
Let $D$ be a link diagram, and $R_1$, $R_2$ regions of $D$ $(R_1 \ne R_2)$. 
Let $c$ be a crossing point of $D$. 
If $c$ satisfies $c\in \partial R_1 \cap \partial R_2$, then 
the region crossing changes on $R_1$ and $R_2$ do not change $c$. 
\label{r1r2}
\end{lem}
\phantom{x}

\begin{figure}[ht]
\begin{center}
\includegraphics[width=80mm]{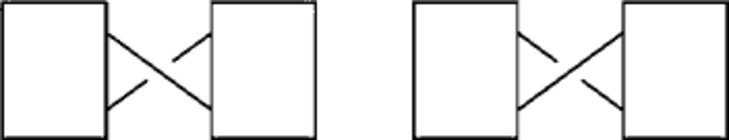}
\caption{}
\label{reduced}
\end{center}
\end{figure}

\begin{figure}[ht]
\begin{center}
\includegraphics[width=25mm]{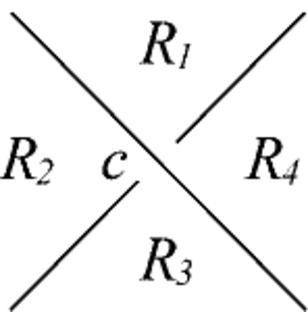}
\caption{}
\label{4r}
\end{center}
\end{figure}

\noindent A link diagram $D$ on $S^2$ is \textit{reducible} if $D$ has a crossing as shown in Figure \ref{reduced}, 
where each square means a diagram of a tangle. 
A link diagram $D$ on $S^2$ is \textit{reduced} if $D$ is not reducible. 
We call such a crossing a \textit{reducible crossing}, 
and the set of a reducible crossing and one of the squares a \textit{reducible part}. 
We have the following proposition: 

\phantom{x}
\begin{prop}
A link diagram $D$ is a reducible link diagram if and only if  there exists a crossing $c$ of $D$ such that 
the regions $R_1,R_2,R_3$ and $R_4$ around $c$ as shown in Figure \ref{4r} satisfy $R_1=R_3$ or $R_2=R_4$. 
\label{non-reduced}
\end{prop}
\phantom{x}

\noindent We can shade some regions of $D$ so that each two regions which are adjacent by an edge of $\vert D \vert$ are shaded and unshaded. 
We call such shading a checkerboard coloring. 
From Lemma \ref{r1r2} and Proposition \ref{non-reduced}, we have the following corollary:

\phantom{x}
\begin{cor}
Let $D$ be a reduced link diagram with a checkerboard coloring, 
and $D'$ the diagram obtained from $D$ by region crossing changes at all the shaded regions. 
Then, $D=D'$. 
\label{checker}
\end{cor}
\phantom{x}

\noindent We remark that Corollary \ref{checker} does not hold for a reducible link diagram 
(see, for example, the diagram $E$ in Figure \ref{rcc-ex}). 
From Corollary \ref{checker}, we have the following corollary: 

\phantom{x}
\begin{cor}
Let $D$ be a reduced link diagram, 
and $B$ the set of all the regions of $D$ shaded in a checkerboard coloring. 
Let $P$ be a subset of $B$ consisting of non-empty regions of $D$. 
Then, $D(P)=D(B-P)$. 
\label{ho-checker}
\end{cor}
\phantom{x}

\noindent From Corollary \ref{ho-checker}, we have the following corollary:

\phantom{x}
\begin{cor}
Let $D$ be a reduced link diagram, and $P$ a set of regions of $D$. 
Then there exist just one or three sets $P_i$ $(i=1$ or $i=1,2,3)$ of regions of $D$ such that $D(P)=D(P_i)$, 
where $P\ne P_i$, $P_i\ne P_j$ $(i\ne j$, $i,j=1,2,3)$. 
\end{cor}
\phantom{x}

\noindent From Lemma \ref{r1r2} and that $D(R,R)=D$, we have the following corollary:

\phantom{x}
\begin{cor}
Let $D$ be a link diagram, and $c$ a crossing point of $D$. 
Let $R_1,R_2,R_3$ and $R_4$ be regions of $D$ around $c$  as shown in Figure \ref{4r}. 
If $R_1\ne R_3$ and $R_2\ne R_4$, then the region crossing changes at $R_1, R_2, R_3$ and $R_4$ do not change $c$. 
\label{4cor}
\end{cor}
\phantom{x}

\section{Proof of Theorem \ref{mainthm}}
\label{proof-s}

In this section, we prove Theorem \ref{mainthm}.

\phantom{x} 
\noindent {\it Proof of Theorem \ref{mainthm}.}
Let $D$ be a knot diagram, and $c$ a crossing point of $D$. 
We show that we can make the crossing change at $c$ by region crossing changes by an induction on the number $k$ of reducible crossings of $D$. 
If $k=0$, i.e., $D$ is a reduced diagram, we can obtain the regions of $D$ such that we can change only $c$ by region crossing changes at the regions by the following procedure: \\

\begin{figure}[ht]
\begin{center}
\includegraphics[width=45mm]{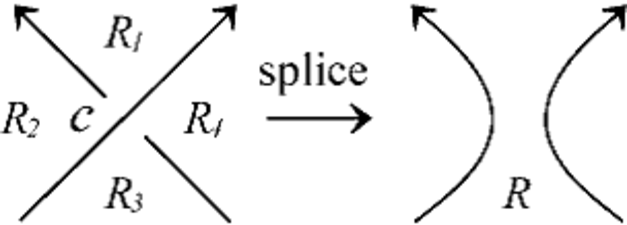}
\caption{}
\label{splice}
\end{center}
\end{figure}

\noindent {\bf Step 1.} We splice $D$ at $c$ by giving $D$ an orientation (see Figure \ref{splice}). 
Then, we obtain a diagram $D_s =D^1 \cup D^2$ of a two-component link. \\

\noindent {\bf Step 2.} We apply a checkerboard coloring for one component $D^1$ of $D_s$ by ignoring another component $D^2$ so that the region $R$ in Figure \ref{splice} is unshaded. \\

\noindent {\bf Step 3.} We take the regions of $D$ corresponding to the shaded regions of $D_s$. 
\begin{figure}[ht]
\begin{center}
\includegraphics[width=90mm]{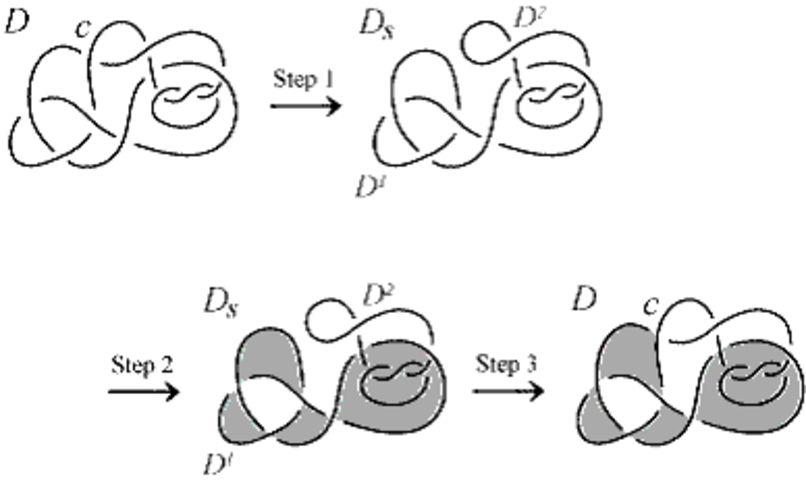}
\caption{}
\label{pf-ex}
\end{center}
\end{figure}
An example of the above procedure is shown in Figure \ref{pf-ex}. 
By Lemma \ref{r1r2} and Corollary \ref{4cor}, a crossing point of $D$ which corresponds to 
a self-crossing point of $D^1$ or $D^2$ is not changed by the region crossing changes at the regions. 
By Lemma \ref{r1r2}, a crossing point of $D$ which corresponds to a crossing point between $D^1$ and $D^2$ is 
not changed by the region crossing changes at the regions. 
Hence we can change only $c$. 
Therefore, the theorem holds for reduced knot diagrams. \\

We remark that if $D_s$ has a reducible crossing $d$ $(\ne c)$, $d$ corresponds to a crossing on $\partial R_1 \cap \partial R_3$ in $D$. 
That is why we apply the checkerboard coloring so that $R$ is unshaded in Step 2. \\

\begin{figure}[ht]
\begin{center}
\includegraphics[width=45mm]{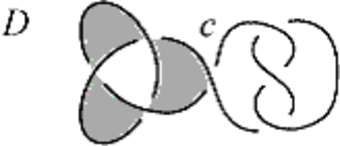}
\caption{}
\label{red-c-case}
\end{center}
\end{figure}

Here, we consider a special case $D$ has just one reducible crossing and $c$ is the reducible crossing. 
Apply the checkerboard coloring to one reducible part as shown in Figure \ref{red-c-case}. 
Then we can change only $c$ by the region crossing changes at the shaded regions. 

\begin{figure}[ht]
\begin{center}
\includegraphics[width=70mm]{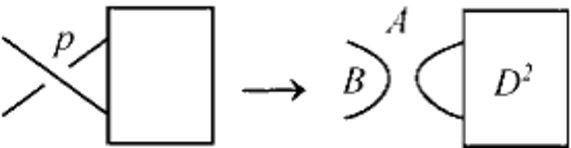}
\caption{}
\label{p-ab}
\end{center}
\end{figure}

We next consider the other case. 
We assume that the theorem holds for all the knot diagrams with $k$ reducible crossings ($k \ge 0$). 
Now we consider knot diagrams with $k+1$ reducible crossings. 
In this case, there exists a reducible crossing $p$ of an innermost reducible part $S$ which does not include $c$. 
By splicing $D$ at $p$, we obtain a non-connected link diagram consisting of a knot diagram $D^1$ with $c$ and $k$ reducible crossings and a reduced knot diagram $D^2$. 
By the assumption, $D^1$ has regions such that we can change only $c$ by the region crossing changes at the regions. 
We call such set of regions $P$. 
We obtain a set $Q$ of regions of $D$ from $P$ by the following rules: 
Let $A$ be the region of $D^1$ which includes $D^2$, and $B$ the opposite region of $D^1$ (see Figure \ref{p-ab}). 
Let $Q$ includes regions corresponding to $P\setminus {A}$, and \\

\begin{figure}[ht]
\begin{center}
\includegraphics[width=120mm]{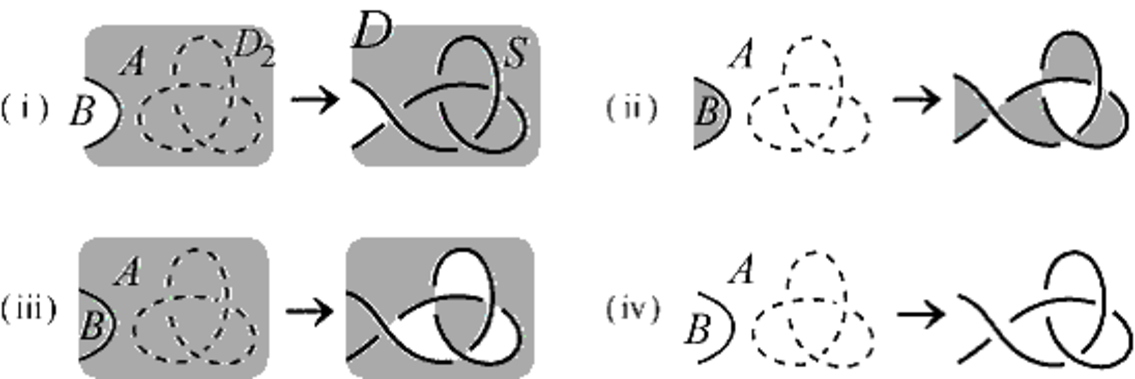}
\caption{}
\label{ab-four}
\end{center}
\end{figure}

\noindent (i) If $A\in P$ and $B\not\in P$, all the regions of $S$ in $A$ is in $Q$. \\
\noindent (ii) If $A\not\in P$ and $B \in P$, the shaded regions of $S$ with the checkerboard coloring such that the outer region of $S$ is white are in $Q$. \\
\noindent (iii) If $A \in P$ and $B \in P$, then the shaded regions of $S$ in $A$ with the checkerboard coloring such that the outer region of $S$ is black are in $Q$. \\
\noindent (iv) If $A \not\in P$ and $B \not\in P$, all the regions of $S$ in $A$ is not in $Q$ (see Figure \ref{ab-four}). \\
Then, $Q$ is the set of regions which change only $c$. 
\hfill$\square$
\phantom{x}

\section{Region unknotting number}
\label{run}

In this section, we discuss the region unknotting number of a knot diagram and a knot. 
We have the following lemma:

\phantom{x}
\begin{lem}
Let $D$ be a reduced knot diagram, and $c(D)$ the crossing number of $D$. 
Then we have
$$u_R(D)\le \frac{c(D)}{2}+1.$$

\phantom{x}
\begin{proof}
For a reduced knot diagram $D$ with a checkerboard coloring, we denote by $b$ (resp. $w$) the number of regions colored black (resp. white). 
We have 
\begin{align*}
u_R(D)&\le \lfloor \frac{b}{2}\rfloor +\lfloor \frac{w}{2}\rfloor \\
&\le \frac{b+w}{2}
\end{align*}
because of Corollary \ref{ho-checker}, 
where $\lfloor x\rfloor =\max \{ n\in \mathbb{Z} \vert n\le x\} $. 
Since $b+w$ means the number of regions of $D$, we have 
$$u_R(D)\le \frac{c(D)+2}{2}.$$
\end{proof}
\label{urd-lem}
\end{lem}
\phantom{x}

\begin{rem}
From the proof of Lemma \ref{urd-lem}, it can also be said that the region unknotting number of a reduced knot diagram $D$ is 
less than or equal to half the number of regions of $D$. 
\end{rem}
\phantom{x}

\begin{rem}
The equality in Lemma \ref{urd-lem} does not hold if $c(D)$ is even and both $b$ and $w$ are odd, or $c(D)$ is odd. 
\end{rem}
\phantom{x}

\noindent We show an example of region unknotting numbers of knot diagrams.

\begin{figure}[ht]
\begin{center}
\includegraphics[width=130mm]{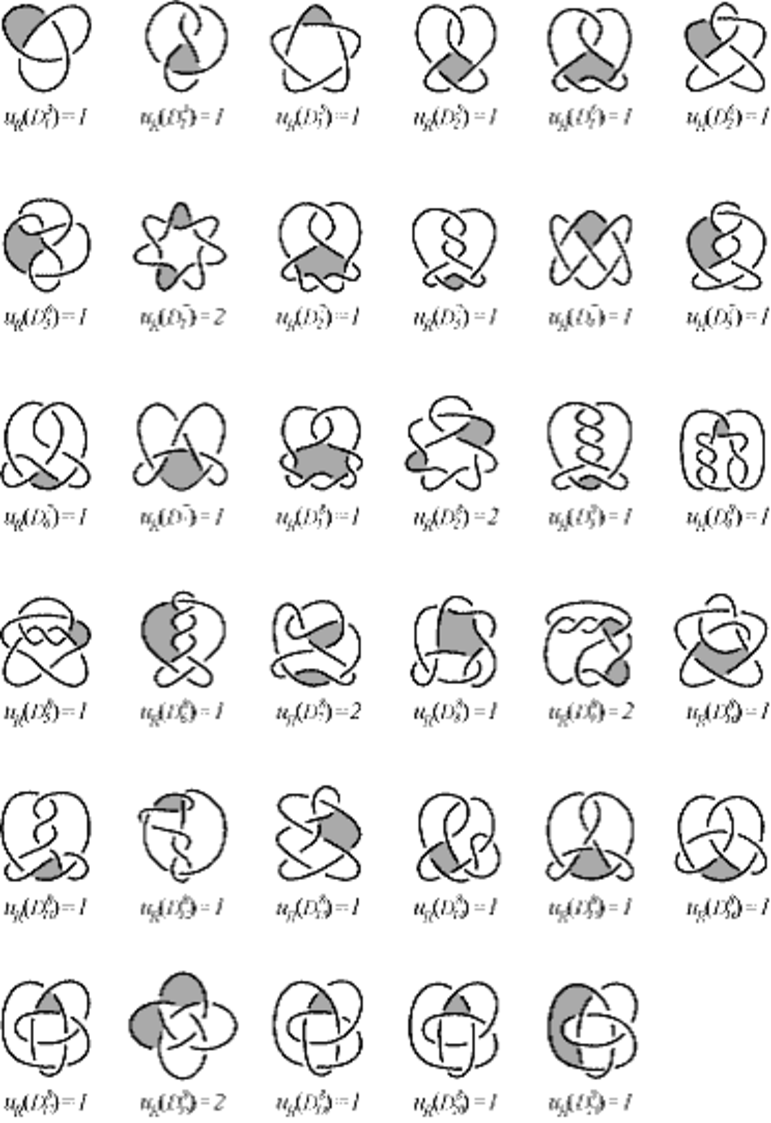}
\caption{}
\label{table}
\end{center}
\end{figure}

\phantom{x}
\begin{eg}
In Figure \ref{table}, we list all the knot diagrams based on Rolfsen's knot table \cite{rolfsen} with crossing number eight or less 
and their region unknotting numbers. 
We denote by $D^m_n$ the diagram of $m_n$ in Rolfsen's knot table (for example, we denote by $D^3_1$ the diagram of $3_1$). 
\end{eg}
\phantom{x}

\noindent We prove Theorem \ref{urk-thm} by using Lemma \ref{urd-lem}:\\

\phantom{x}
\noindent \textit{Proof of Theorem \ref{urk-thm}}. 
For a knot $K$ and a minimal crossing diagram $D$ of $K$, we have 
$$u_R(K)\le u_R(D)\le \frac{c(D)}{2}+1=\frac{c(K)}{2}+1$$
because $D$ is a reduced knot diagram. 
\hfill$\square$\\
\phantom{x}

\noindent In the following example, we show the region unknotting numbers of all the prime knots with crossing number nine or less: 

\phantom{x}
\begin{eg}
The knots $7_1,8_2,8_7,8_9,8_{18},9_1,9_3,9_6, 9_{35},9_{40}$ have region unknotting numbers two. 
The other prime knots with crossing number nine or less have the region unknotting number one. 
\label{run-9ex}
\end{eg}
\phantom{x}

\begin{rem}
For the above knots in Example \ref{run-9ex}, the region unknotting numbers are realized by the diagrams in Rolfsen's knot table. 
We note that the knots $7_1,8_2,8_7,8_9,8_{18},9_1,9_3,9_{35},9_{40}$ have only one minimal crossing diagrams, respectively 
up to horizontal mirror image and vertical mirror image, 
where the horizontal mirror image of a knot diagram $D$ is obtained from $D$ by reflecting $D$ across a vertical plane, 
and the vertical mirror image of $D$ is obtained from $D$ by changing all the crossings of $D$. 
The knot $9_6$ has just two minimal crossing diagrams, whose region unknotting numbers are two (see Figure \ref{96-2}). 
We will discuss how to obtain all minimal crossing diagrams on $S^2$ of prime alternating knots in the appendix.
\begin{figure}[ht]
\begin{center}
\includegraphics[width=70mm]{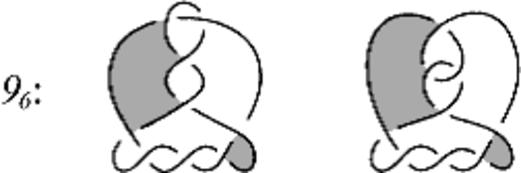}
\caption{}
\label{96-2}
\end{center}
\end{figure}
We remark that there exist minimal crossing diagrams of prime alternating knots which do not realize the region unknotting numbers. 
For example, J. Banks suggested that the two minimal crossing diagrams $D$ and $E$ of $9_{26}$ in Figure \ref{jb-ex} have 
$u_R(D)=1$ and $u_R(E)=2$. 
\begin{figure}[ht]
\begin{center}
\includegraphics[width=70mm]{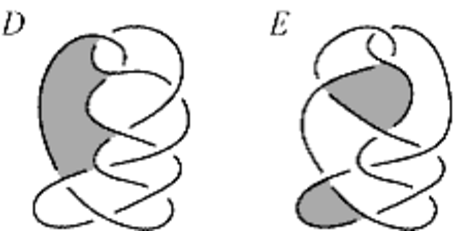}
\caption{}
\label{jb-ex}
\end{center}
\end{figure}
\label{jb-rem}
\end{rem}
\phantom{x}

\noindent For twist knots, we have the following proposition: 

\phantom{X}
\begin{prop}
A twist knot $K$ has $u_R(K)=1$. 

\phantom{x}
\begin{proof}
From the minimal crossing diagram of $K$ in Figure \ref{twist}, 
we can obtain a diagram of the trivial knot by a region crossing change at the region $P$ or $Q$. 
\begin{figure}[ht]
\begin{center}
\includegraphics[width=30mm]{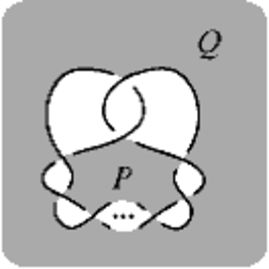}
\caption{}
\label{twist}
\end{center}
\end{figure}
\end{proof}
\end{prop}
\phantom{x}

\noindent For $(2,2n+1)$-torus knots, we have the following proposition:

\phantom{x}
\begin{prop}
If a knot $K$ is the $(2,4m-1)$-torus knot or the $(2,4m+1)$-torus knot $(m=1,2,\dots )$, then $u_R(K)=m$. 

\phantom{x}
\begin{proof}
First of all, we remark that the $(2,2n+1)$-torus knot has only one minimal crossing diagram on $S^2$ 
as shown in Figure \ref{n-foil} (see the appendix), 
and the region crossing change at $P$ or $Q$ in Figure \ref{n-foil} is of no use because it always transforms a diagram 
into just the vertical mirror image. 
\begin{figure}[ht]
\begin{center}
\includegraphics[width=30mm]{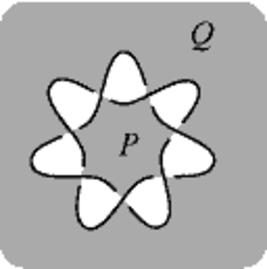}
\caption{}
\label{n-foil}
\end{center}
\end{figure}
We prove that the $(2,4m-1)$-torus knot has the region unknotting number $m$ by an induction. 
When $m=1$, the trefoil knot has the region unknotting number one. 
We assume that in the case of $m=k$, the $(2,4k-1)$-torus knot $K$ has the region unknotting number $k$. 
We shall prove for the case of $m=k+1$, that the $(2,4k+3)$-torus knot $K'$ has the region unknotting number $k+1$. \\
$\bullet$ First, we prove that $u_R(K')\le k+1$. 
Let $D$ be the minimal crossing diagram of $K$. 
We add two full twists to a pair of edges of $D$ which bound a bigonal region as shown in Figure \ref{replace} 
so that we obtain the minimal crossing diagram $D'$ of $K'$. 
\begin{figure}[ht]
\begin{center}
\includegraphics[width=35mm]{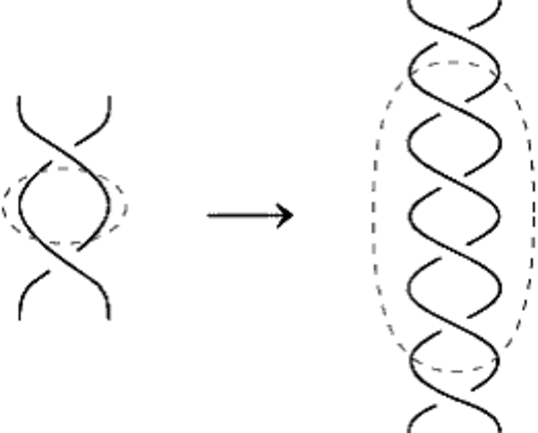}
\caption{}
\label{replace}
\end{center}
\end{figure}
Since we obtain $D$ from $D'$ by a region crossing change at any region in the two full twists and 
Reidemeister moves of type I\hspace{-0.5pt}I (see Figure \ref{re-replace}), 
we have $u_R(K')\le u_R(K)+1=k+1$. \\
\begin{figure}[ht!]
\begin{center}
\includegraphics[width=50mm]{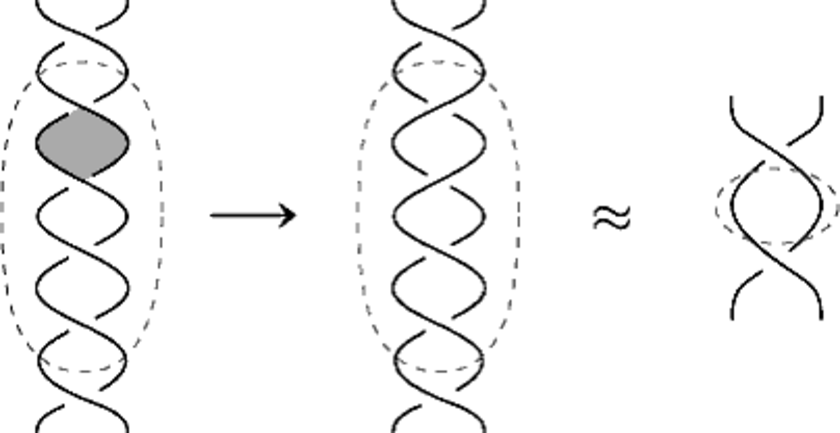}
\caption{}
\label{re-replace}
\end{center}
\end{figure}
$\bullet$ Next, we prove $u_R(K')\ge k+1$ by an indirect proof. 
We assume that $u_R(K')\le k$. 
Let $D'$ be the minimal crossing diagram of $K'$. 
Let $R'$ be a set of $u_R(K')$ regions of $D'$ such that $D'(R')$ represents the trivial knot. 
Since the number of bigonal regions of $D'$ is $4k+3$ and $u_R(K')=u_R(D')\le k$, there exist connected four bigonal regions 
satisfying the following condition: 
One region ${R_1}'$ of them is in $R'$, and the other three regions of them are not in $R'$. 
By applying a region crossing change at ${R_1}'$, we obtain from $D'$ a diagram $D'({R_1}')$ which represents the knot $K$, 
and by applying Reidemeister moves of type I\hspace{-0.5pt}I, we obtain from $D'({R_1}')$ a minimal crossing diagram $D$ of $K$. 
By region crossing changes at the regions of $D$ which corresponds to the regions $R'-{R_1}'$ of $D'$, 
we obtain from $D$ a diagram representing the trivial knot. 
Hence $u_R(K)\le u_R(K')-1\le k-1$ which contradicts $u_R(K) = k$. 
Hence we have $u_R(K') \ge k+1$, and therefore $u_R(K') = k+1$. 
The $(2,5)$-torus knot has the region unknotting number one, and 
we can prove similarly for $(2,4m+1)$-torus knots. 
\end{proof}
\label{n-foil-prop}
\end{prop}
\phantom{x}

\noindent From Proposition \ref{n-foil-prop}, we have the following corollary: 

\phantom{x}
\begin{cor}
For an arbitrary non-negative integer $n$, there exists a knot $K$ which satisfies $u_R(K)=n$. 
\end{cor}

\section*{Appendix}

\begin{figure}[ht!]
\begin{center}
\includegraphics[width=75mm]{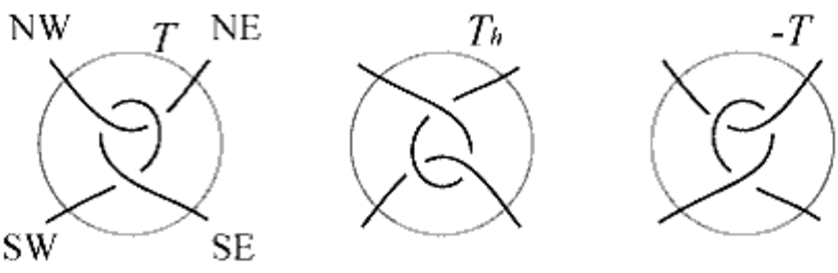}
\caption{}
\label{tangle-def}
\end{center}
\end{figure}

\begin{figure}[ht!]
\begin{center}
\includegraphics[width=45mm]{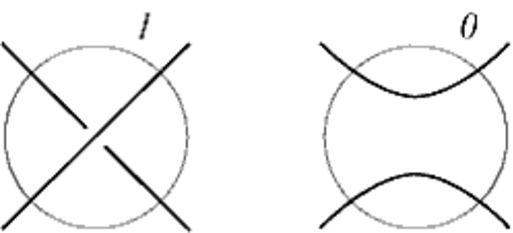}
\caption{}
\label{one-tangle}
\end{center}
\end{figure}

\begin{figure}[ht!]
\begin{center}
\includegraphics[width=80mm]{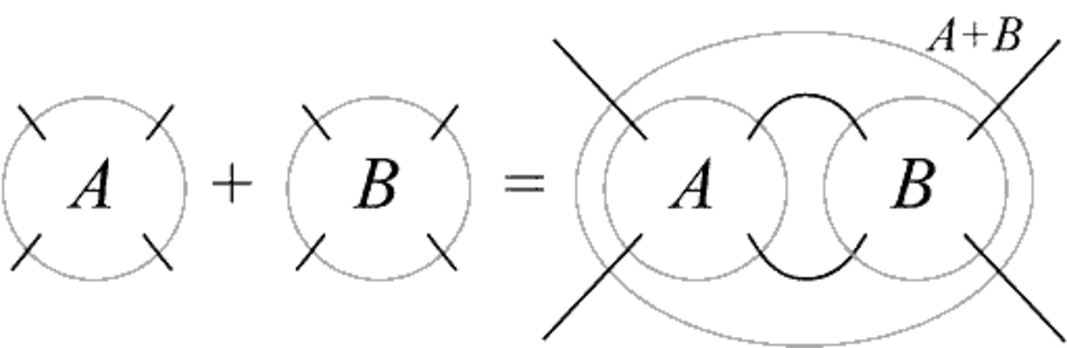}
\caption{}
\label{sum-def}
\end{center}
\end{figure}

In this appendix, we explain how to obtain another minimal crossing diagram by a flyping 
from a minimal crossing diagram of a prime alternating knot. 
Then we show how to obtain all the minimal crossing diagrams of a prime alternating knot. 
In this appendix, a \textit{tangle} is a portion of a knot diagram 
from which there emerge just four arcs pointing in the four compass directions NW, NE, SW, and SE. 
For a tangle $T$, we denote by $T_h$ (resp. $T_v$) the result of rotation in a horizontal (resp. vertical) axis, 
and $-T$ that of crossing changes at all the crossing points of $T$ as shown in Figure \ref{tangle-def}. 
We denote by $1$ (resp. $0$) the tangle with one crossing (resp. no crossings) 
as shown in the left side (resp. right side) of Figure \ref{one-tangle}. 
For two tangles $A$ and $B$, we define the \textit{sum $A+B$ of $A$ and $B$} to be the result of the operation of Figure \ref{sum-def}. 
We also denote by $A-B$ the sum of $A$ and $-B$. 
A tangle $T$ is a \textit{tangle sum} if $T=T_1+T_2$, where neither $T_1$ nor $T_2$ is the tangle $0$. 

\begin{figure}[ht!]
\begin{center}
\includegraphics[width=130mm]{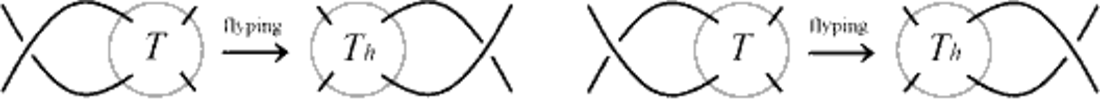}
\caption{}
\label{flyping-def}
\end{center}
\end{figure}

\begin{figure}[ht!]
\begin{center}
\includegraphics[width=75mm]{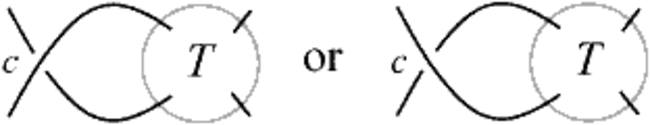}
\caption{}
\label{c-t-take}
\end{center}
\end{figure}

\begin{figure}[ht!]
\begin{center}
\includegraphics[width=80mm]{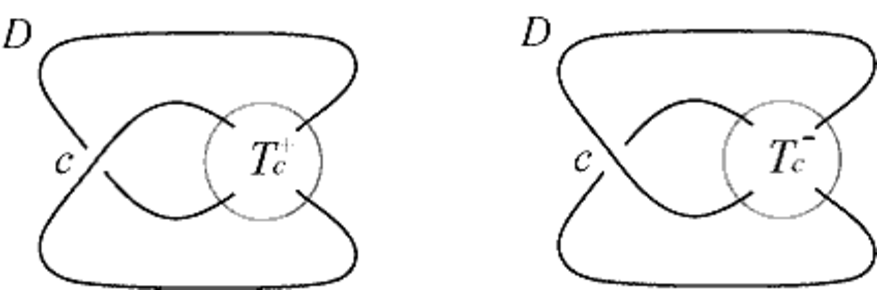}
\caption{}
\label{t-c-tangle}
\end{center}
\end{figure}

\begin{figure}[ht!]
\begin{center}
\includegraphics[width=90mm]{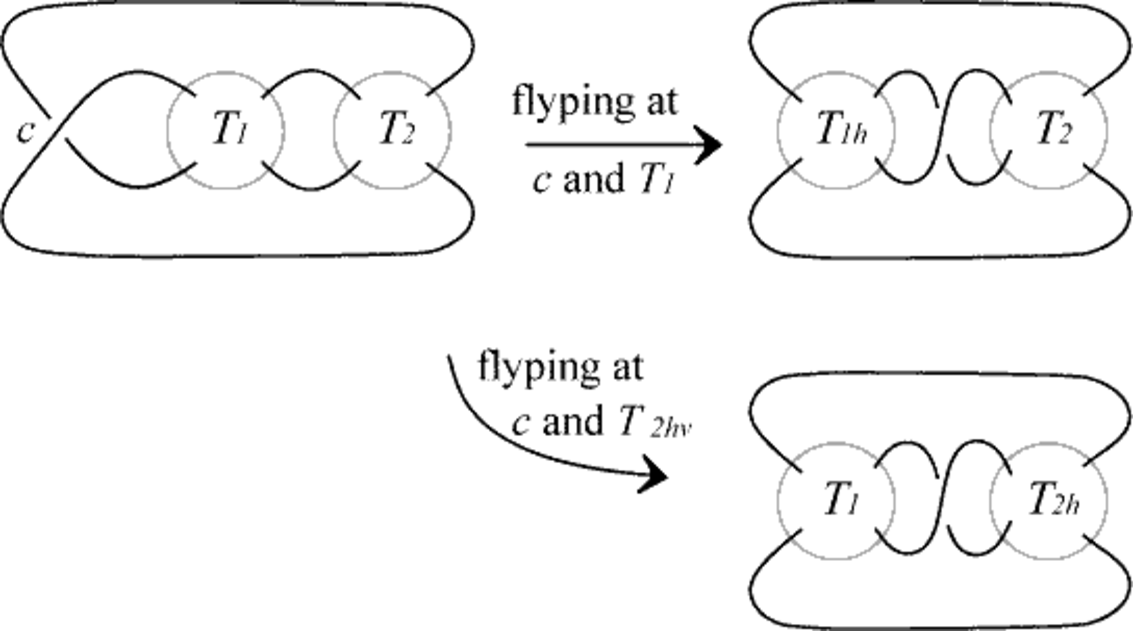}
\caption{}
\label{t1-t2-v-m}
\end{center}
\end{figure}

\begin{figure}[ht!]
\begin{center}
\includegraphics[width=85mm]{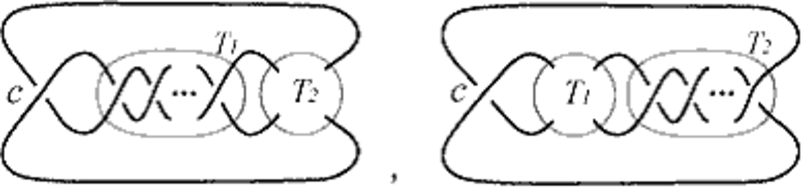}
\caption{}
\label{half-twist}
\end{center}
\end{figure}

\begin{figure}[ht!]
\begin{center}
\includegraphics[width=125mm]{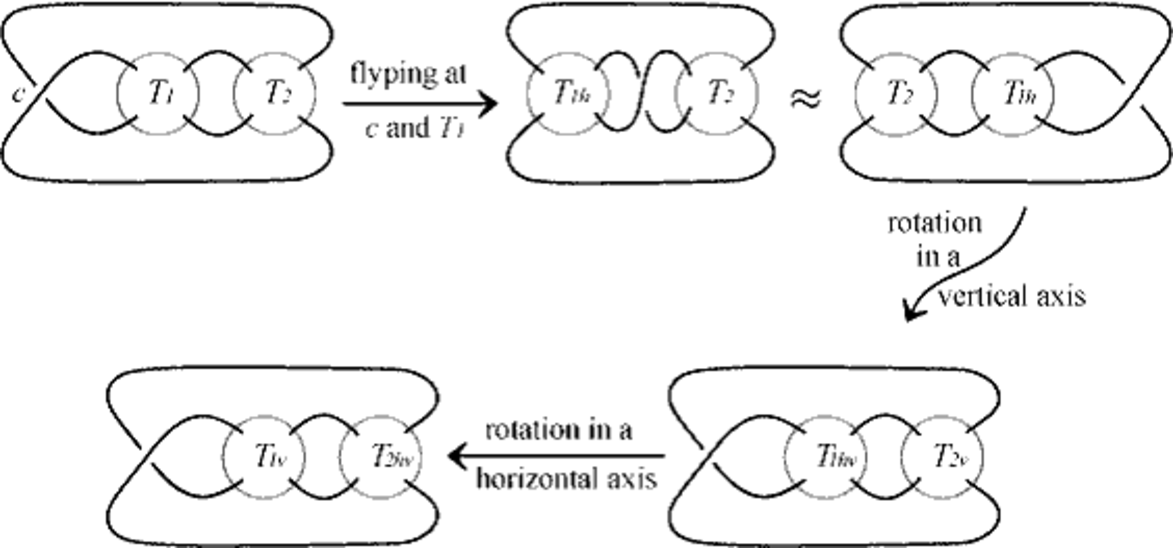}
\caption{}
\label{hv-eq}
\end{center}
\end{figure}

\noindent Let $D$ be a knot diagram which includes a tangle $1+T$ or $-1+T$. 
\textit{Flyping} is a local transformation on $D$ which replaces $1+T$ by $T_h+1$, or $-1+T$ by $T_h-1$ as shown in Figure \ref{flyping-def}. 
W. Menasco and M. Thistlethwaite showed that Tait's third conjecture is true, that is, we can change $D_1$ into $D_2$ 
by performing a finite number of flypings for any two reduced alternating diagrams $D_1$ and $D_2$ of an alternating knot $K$ \cite{me-thi}. 
Hence we can obtain all the minimal crossing diagrams of a prime alternating knot $K$ from a minimal crossing diagram of $K$ by flypings. 
Let $D$ be a minimal crossing diagram of a non-trivial knot, 
and $c$ a crossing point of $D$. 
Let $T$ be a tangle in $D$ whose NW arc and SW arc meet at $c$ as shown in Figure \ref{c-t-take}. 
Since $c$ can be considered as the tangle $1$ or $-1$, 
we can apply a flyping there. 
We call such a flyping \textit{flyping at $c$ and $T$}. 
A flyping on a knot diagram is \textit{trivial} if we obtain $D$ or the (vertical, horizontal, or vertical and horizontal) mirror image of $D$ 
from $D$ by the flyping. 
When we can apply a non-trivial flyping at a crossing point $c$ and a tangle $T$ of a diagram $D$, we say that $D$ \textit{admits} 
a non-trivial flyping at $c$ (and $T$). 
Now we explain how to obtain all the tangles $T$ such that $D$ admits non-trivial flypings at a crossing point $c$ and $T$ for a knot diagram $D$. \\
\noindent From a knot diagram $D$ on $S^2$ and a crossing point $c$ of $D$, we obtain two tangles $T^+_c$ and $T^-_c$ 
such that we obtain $D$ from $1+T^+_c$ (resp. $-1+T^-_c$) by connecting the NW arc and the NE arc, and the SW arc and the SW arc 
$($see Figure \ref{t-c-tangle}$)$, where we remark that the tangle $1$ (resp. $-1$) corresponds to $c$. 
We note that $T^+_c$ and $T^-_c$ are unique by regarding that 
$\pm 1 +{T^{\varepsilon}_c}_{hv}$ is equivalent to $\pm 1+T^{\varepsilon}_c$ on $S^2$ $(\varepsilon =+,-)$. 
We remark that a flyping at $c$ and $T^+_c$ or $T^-_c$ is trivial because it comes out just horizontal mirror image of $D$. 
Hence if neither $T^+_c$ nor $T^-_c$ is a tangle sum, the diagram $D$ does not admit a non-trivial flyping at $c$. 
When $T^{\varepsilon}_c$ is a tangle sum of tangles $T_1$ and $T_2$ $(\varepsilon =+,-)$, we can apply a flyping at $c$ and $T_1$. 
Remark that a flyping at $c$ and $T_1$ is equivalent to a flyping at $c$ and ${T_2}_{hv}$ up to horizontal mirror image 
$($see Figure \ref{t1-t2-v-m}$)$. 
Then we shall consider flypings at $c$ and only $T_1$. 
For $T^{\varepsilon}_c=T_1+T_2$, the flyping at $c$ and $T_1$ is trivial if $T_1$ or $T_2$ is the sum of some tangles 
$\varepsilon 1$ $(\varepsilon =+,-)$ as shown in Figure \ref{half-twist}. 
The flyping at $c$ and $T_1$ is also trivial if $T_1$ and $T_2$ satisfy $T_{1hv}=T_1$ and $T_{2v}=T_2$, 
or $T_{1v}=T_1$ and $T_{2hv}=T_2$ (see Figure \ref{hv-eq}). 
Then, to obtain all the minimal crossing diagrams on $S^2$ of a reduced alternating knot $K$, 
we will consider all the flypings at all the crossing points $c$ and all the tangles $T_1$ such that $T^{\varepsilon}_c=T_1+T_2$ $(\varepsilon =+,-)$ 
of a minimal crossing diagram $D$ of $K$ except the following three cases: \\

\begin{description}
\item[$($i$)$] the tangle $T^{\varepsilon}_c$ is not a tangle sum $(\varepsilon =+,-)$, 
\item[$($ii$)$] the tangle $T_1$ or $T_2$ is the sum of some tangles $\varepsilon 1$ $(\varepsilon =+,-)$, 
\item[$($iii$)$] the tangles $T_1$ and $T_2$ satisfy $T_{1hv}=T_1$ and $T_{2v}=T_2$, or $T_{1v}=T_1$ and $T_{2hv}=T_2$. 
\end{description}

\noindent In the following examples, we will find all the minimal crossing diagrams of some knots by the above procedure. 

\phantom{x}
\begin{eg}
A $(2,2n+1)$-torus knot has only one minimal crossing diagram $D$ on $S^2$: 
Let $D$ in Figure \ref{ex-5-tc} be the minimal crossing diagram of a $(2,2n+1)$-torus knot $(n=1,2,\dots )$. 
For every crossing point $c$, we obtain the same tangles $T^+_c$ and $T^-_c$ as shown in Figure \ref{ex-5-tc}. 
\begin{figure}[ht!]
\begin{center}
\includegraphics[width=100mm]{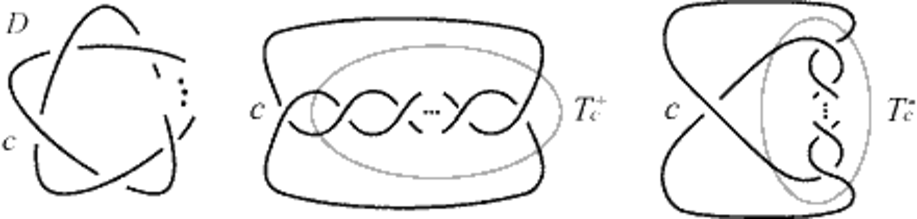}
\caption{}
\label{ex-5-tc}
\end{center}
\end{figure}
The tangle $T^+_c$ is a sum of two tangles $T_1$ and $T_2$, where $T_1$ is $k$ half twists and $T_2$ is $2n-k$ half twists $(k=1,2,\dots , 2n-1)$. 
The tangle $T^-_c$ is not a tangle sum. 
Then $T^+_c$ and $T^-_c$ satisfy the cases (ii) and (i), respectively. 
therefore, we can not apply non-trivial flypings on $D$. 
\end{eg}
\phantom{x}

\noindent Similarly, we have the following example: 

\begin{eg}
A knot $K$ with Conway's notation $mn$ or $m,n$ has only one minimal crossing diagram on $S^2$ 
$(m,n\ne 0 \in \mathbb{Z}, mn>0)$. 
\end{eg}
\phantom{x}

\noindent We show that the knot $8_2$ has only one minimal crossing diagram on $S^2$: 

\phantom{x}
\begin{eg}
For the minimal crossing diagram $D$ of the knot $8_2$ in Figure \ref{d-82}, we call each crossing point $a,b,\dots $ and $h$ as shown in the figure. 
\begin{figure}[ht!]
\begin{center}
\includegraphics[width=30mm]{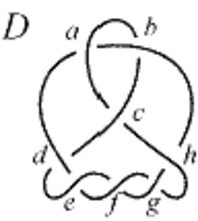}
\caption{}
\label{d-82}
\end{center}
\end{figure}
Neither $T_a^-$, $T_b^-$, $T_c^+$, $T_d^-$, $T_e^-$, $T_f^-$, $T_g^-$ nor $T_h^-$ is a tangle sum, 
$T_a^+$ and $T_b^+$ are the tangle sums which consist of a tangle and the tangle $1$, 
and $T_c^-$ is the tangle sum of the two tangles with vertical twists $T_1$ and $T_2$ satisfying $T_{1hv}=T_1$ and $T_{2v}=T_2$. 
The tangles $T_d^+$, $T_e^+$, $T_f^+$, $T_g^+$ and $T_h^+$ are always tangle sums $T_1+T_2$ 
such that $T_1$ or $T_2$ is the sum of $k$ tangles $1$ $(k\ge 1)$ (see Figure \ref{82-12}). 
Hence $D$ admits no non-trivial flypings, and therefore $8_2$ has only one minimal crossing diagram $D$ on $S^2$. 
\begin{figure}[ht!]
\begin{center}
\includegraphics[width=125mm]{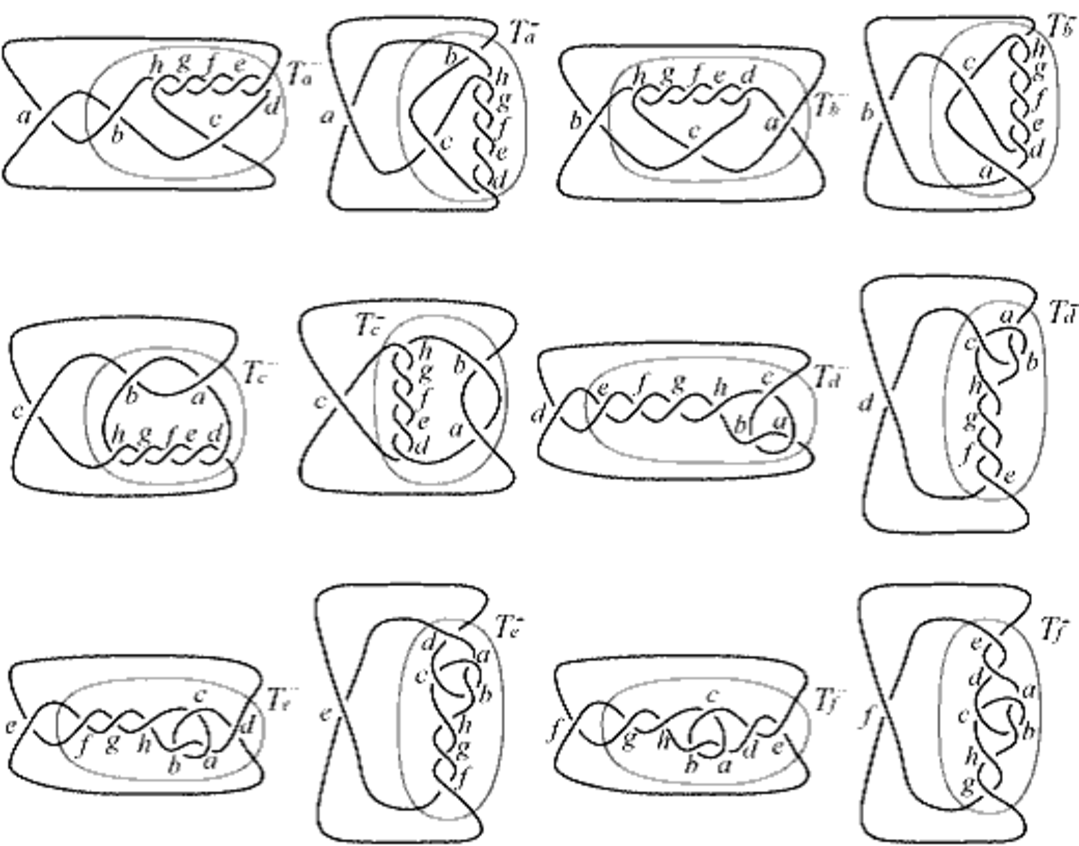}
\caption{}
\label{82-12}
\end{center}
\end{figure}
\end{eg}
\phantom{x}

\noindent We next show that the knot $9_6$ has just two minimal crossing diagrams on $S^2$: 

\phantom{x}
\begin{eg}
For the minimal crossing diagram $D$ of $9_6$ in Figure \ref{96-cd-1}, we call each crossing $a,b,\dots$ and $i$ as shown in the figure. 
Then $D$ admits non-trivial flypings only at $c$ and the tangle $T^c_1$ and at $d$ and $T^d_1$, 
and we obtain another diagram $D'$ in Figure \ref{96-cd-2} by the flyping at $c$ and $T^c_1$ and at $d$ and $T^d_1$, respectively, 
where we denote by $T^-_c=T^c_1+T^c_2$ and $T^-_d=T^d_1+T^d_2$ the tangles in Figure \ref{96-cd-1}. 
\begin{figure}[ht!]
\begin{center}
\includegraphics[width=105mm]{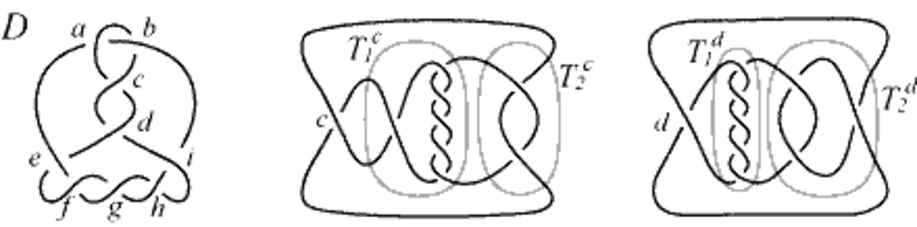}
\caption{}
\label{96-cd-1}
\end{center}
\end{figure}
The diagram $D'$ with the crossing points $a',b',\dots$ and $i'$ in Figure \ref{96-cd-2} admits non-trivial flypings only at 
$c'$ and ${T'}^{c'}_{11}$, $c'$ and ${T'}^{c'}_{12}$, $d'$ and ${T'}^{d'}_{11}$, and $d'$ and ${T'}^{d'}_{12}$ as depicted in Figure \ref{96-cd-2}, 
and we obtain the diagram $D$ in Figure \ref{96-cd-1} by flyping there, respectively. 
Hence $9_6$ has just two minimal crossing diagrams $D$ and $D'$ on $S^2$. 
\begin{figure}[ht!]
\begin{center}
\includegraphics[width=110mm]{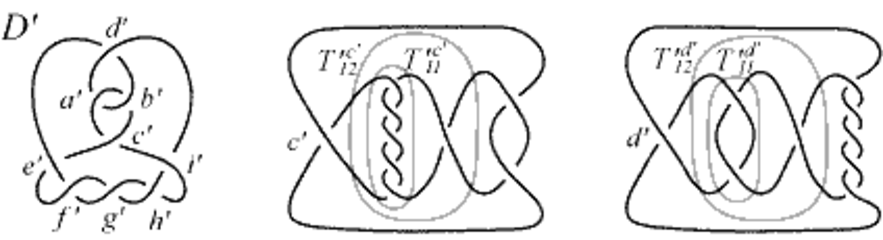}
\caption{}
\label{96-cd-2}
\end{center}
\end{figure}
\end{eg}

\section*{Acknowledgments}
The author is grateful to Akio Kawauchi, Taizo Kanenobu, and the members of Friday Seminar on Knot Theory 
in Osaka City University who assisted her in helpful advice, valuable discussions and the tender encouragement. 
She especially thanks Kengo Kishimoto for giving her such an interesting question and valuable advice and discussions. 
She also thanks Kenneth C. Millett and Makoto Ozawa for helpful comments, advice and information. 
She is deeply grateful to Jessica E. Banks for her valuable comments and many helpful suggestions. 
She is also grateful to the referee for careful reading and helpful advice. 
She was supported by JSPS Research Fellowships for Young Scientists. 

\maketitle

\end{document}